

\documentstyle{amsppt} 
\magnification=\magstep1 
\NoRunningHeads

\vsize=7.4in     



\def \qed {\vrule height6pt  width6pt depth0pt} 

\def\X{\frak X}
\def\e{\tau}
\def\d{\delta}
\def\a{\alpha}
\def\Sg{\Sigma}
\def\S+{\Sigma^+}
\def\n{||}


\topmatter

\title 
The Complete Continuity Property and 
Finite Dimensional Decompositions 
\endtitle 
\author 
Maria Girardi
\quad and \quad
William B. Johnson
\endauthor
\address 
University of South Carolina, 
Department of Mathematics, 
Columbia, SC  29208. \linebreak    
E-mail: girardi\@math.scarolina.edu  
\endaddress 
\address 
Texas A\&M University,
Department of Mathematics, 
College Station, TX
77843. \linebreak           
E-mail: WBJ7835\@venus.tamu.edu  
\endaddress 
\thanks
The first author was supported in part by NSF DMS-9204301. 
\endthanks 
\thanks
The second author was supported in part by NSF DMS-9003550. 
\endthanks 
\subjclass
46B22, 46B20, 46B28, 46G99
\endsubjclass 
\abstract  
 A Banach space $\X$ has the complete continuity property  (CCP)
if each  bounded linear operator from $L_1$ into $\X$ 
is completely continuous 
(i.e.,  maps weakly convergent sequences to norm convergent sequences).
The main theorem shows that  
a Banach space failing the CCP  
(resp.,  failing  the CCP and failing cotype)  
has a subspace with a 
finite dimensional decomposition  
(resp.,  basis)    
which fails the CCP. 
\endabstract 

\endtopmatter

\document 

\heading{1. introduction}\endheading \vskip 10pt

Given a property  of Banach spaces which is hereditary, it is natural to ask
whether a Banach space has the property if every subspace with a basis  
(or with a finite dimensional decomposition)  
has the property.  The motivation
for such  questions is of course that it is much easier to deal with Banach
spaces which have a basis (or at least a finite dimensional decomposition) 
than with general spaces.  In this note we consider these questions for
the {\it complete continuity property} (CCP), which means that each 
bounded  linear operator from $L_1$ into the space  is completely
continuous  (i.e., carries weakly convergent sequences into norm convergent
sequences).

The CCP is closely connected with the  Radon-Nikod\'ym property (RNP).   
Since a representable operator is completely continuous, 
the RNP implies the CCP; however,
the Bourgain-Rosenthal space  [BR] has the CCP but not the RNP. 
Bourgain [B1] showed that a space failing the RNP 
has a subspace  with a 
finite dimensional decomposition  which fails the RNP.  
Wessel [W] showed that a space failing the CCP 
has a subspace with a basis which fails the RNP. It is open whether a space
has the RNP (respectively, CCP) if every subspace with a basis has the RNP
(respectively, CCP).

Our main theorem shows that 
if $\X$ fails the CCP, 
then there is an operator $T\: L_1 \to \X$ 
that behaves like the identity operator  $I \: L_1 \to L_1$
on the Haar functions 
$\{ h_j \}$. 
Specifically, 
there is a sequence $\{ x^*_n \}$ in the unit ball of $\X^*$ 
such that  $x^*_n$ keeps the image of 
each Haar functions along the $n^{\text{th}}$-level 
large (i.e.  $x^*_n (T h_{2^n+k}) > \d >0$  )  
and  the natural blocking  
$ \{  \text{sp}( T h_{2^n+k} \: k=1, \ldots, 2^n ) \}_n$ 
of the images of the Haar functions  is a  finite dimensional decomposition   
for some subspace $\X_0$.  
Note that $\X_0$ fails the CCP since $T$ is 
not completely continuous  ($T$ keeps the  
Rachemacher functions larger than $\d$ in norm). 
Thus a space failing the CCP 
has a subspace  with a finite dimensional decomposition which fails the CCP. 
In  the language of Banach space geometry, 
the theorem  says that  
in any Banach space  which fails the CCP
grows  a  
separated $\d$-tree with a difference sequence 
naturally blocking into a finite dimensional decomposition. 
If furthermore the space also fails cotype, 
then modifications produce  a      
separated $\d$-tree growing inside a subspace with a basis. 

Throughout this paper, $\X$ denotes an arbitrary Banach space, 
$\X^*$ the dual space of~$\X$,   
and  $S(\X)$ the unit sphere of $\X$.   
The triple $(\Omega, \Sigma, \mu)$~refers 
to the Lebesgue measure space on  $[0,1]$, 
$\S+$ to the sets in $\Sigma$ with positive measure, 
and $L_1$ to $L_1(\Omega, \Sigma, \mu)$.
All notation and terminology, not otherwise explained, are as in [DU].
  
The authors are grateful to  Michel Talagrand  and Peter Casazza for    
helpful discussions.

\heading{2. operator view-point} \endheading \vskip 10pt

A system 
$\Cal A = \{ A^n_k \in\Sg \: n=0,1,2,\ldots \text{ and } k=1,\ldots,2^n \}$  
is a {\it dyadic splitting} of $A^0_1\in\S+$ if  
each $A^n_k$  is partitioned into the   
two sets $A^{n+1}_{2k-1}$ and $A^{n+1}_{2k}$ 
of equal measure for each admissible $n$ and $k$ .  
Thus  the collection $\pi_n = \{ A^n_k \: k=1, \ldots, 2^n \}$ of sets 
along the $n^{\text{\, th}}$-level  partition $A^0_1$ with 
$\pi_{n+1}$ refining $\pi_n$ and $\mu(A^n_k) = 2^{-n}\mu(A^0_1)$. 
To a dyadic splitting corresponds a (normalized) Haar system 
$\{ h_j \}_{j\ge 1}$ where 
$$
   h_1 = \tfrac{1}{\mu(A^0_1)} \, 1_{A^0_1}
   \text{\qquad and \qquad} 
   h_{2^n+k} = \tfrac{2^n}{\mu(A^0_1)} \, 
               (1_{A^{n+1}_{2k-1}} - 1_{A^{n+1}_{2k}} ) 
$$ 
for $n=0,1,2,\ldots $ and  $k=1,\ldots,2^n $. 

A set $N$ in the  unit sphere of the dual of a 
Banach space $\X$ is said to norm  a subspace $\X_0$ 
within $\e>1$ if for each $x\in\X_0$  there is $x^*\in N$ 
such that $\n x \n \le \e x^* (x)$. 
It is well known and easy to see that  a   
sequence $\{\X_j\}$ of subspaces of $\X$ forms a 
finite dimensional decomposition 
with constant at most $\e$ provided that  for each $n\in\Bbb N$ 
the space  generated by $\{ \X_1, \ldots, \X_n \}$ 
can be normed by a  set from $S(\X_{n+1}^\perp)$ within $\e_n > 1$ 
where $\Pi\e_n\le\e$.  

\proclaim{Theorem 1} 
If the bounded linear operator  $T\: L_1\to\X$ is not completely continuous 
and $\{ \e_n \}_{n\ge 0}$  is a sequence of numbers larger than 1, 
then there exists 
\roster 
\item"{(A)}" 
a dyadic splitting $\Cal A = \{ A^n_k \}$ 
\item"{(B)}"
a sequence $\{ x^*_{t_n}\}_{n\ge0}$ in $S(\X^*)$ 
\item"{(C)}"
a finite set $\{ y^*_{n,i} \}_{i=1}^{p_n}$ in  $S(\X^*)$ for  each $n\ge 0$ 
\endroster 
such that for the Haar system  $\{ h_j \}_{j\ge 1}$   
corresponding to $\Cal A$, for some $\delta>0$, and each $n\ge 0$
\roster 
\item $x^*_{t_n} (T h_{2^n+k})  > \d$ for  $k=1,\dots,2^n$ 
\item $\{ y^*_{n,i} \}_{i=1}^{p_n}$ norms  
      $\text{sp} ( T h_j \:  1 \le j \le 2^n )$ 
       within $\e_n$ 
\item $y^*_{n,i} (T h_{2^n + k})  = 0$ for  
$ k=1,\ldots,2^n$ and $i=1,\dots,p_n$.
\endroster 
Note that if $\Pi \e_n $ is finite, then conditions (2) and (3)  
guarantee that  the natural blocking  
$\{ \text{sp} ( T h_j \:  2^{n-1} < j \le 2^n ) \}_{n\ge 0}$  
forms a finite dimensional decomposition  with constant at most $\Pi \e_n $.
\endproclaim
The proof uses the following standard lemma which, for completeness, we shall
prove later. 
\proclaim{Lemma 2} 
If $A\in\S+$ and $\{ g_i \}_{i=1}^n$ is a finite collection of $L_1$ 
functions, then an extreme point $u$ of the set 
$ C  \ \equiv \  \{ f\in L_1 \: 
  |f| \le 1_A \text{ and }\int_{A} f  g_i \, d\mu = 0
  \text{ for  }i=1, \ldots, n \} $   
has the form $|u| = 1_A $.
\endproclaim 

\demo{Proof of Theorem~1} 
Let $T\: L_1\to\X$ be a norm one  operator  that is not completely continuous. 
Then there is a sequence $\{ r_t \}$ in $L_1$ 
and  a sequence $\{ x^*_t \}$ in $S(\X^*)$ satisfying: 
\roster
\item"{(a)}" $\n r_t \n_{L_\infty} \ \le 1$ 
\item"{(b)}" $r_t$ converges to $0$ weakly in $L_1$ 
\item"{(c)}" $4\d \le x^*_t \, T \, r_t$ for some $\d > 0$ .
\endroster
Consider $T^* x^*_t \in L_\infty$.
Since  $\n r_t(T^* x^*_t) \n_{L_\infty}$ is at most 1, by passing to a
subsequence we may assume that $\{r_t (T^* x^*_t)\}$ converges to some function
$h$ in the  weak-star topology on $L_\infty$.     
Since  $\int h \, d\mu \ge 4 \d$  
the set  $A^0_1 \equiv [h\ge\ 4 \d]$ is in $\S+$. 
(Compare this with  [B2, proposition~5]).

We shall construct, by induction  on the level $n$, 
a dyadic splitting of $A^0_1$ along with 
the desired functional. 
Fix $n\ge0$. 

Suppose we are given a finite dyadic splitting  
$\{ A^m_k \: m=0,\ldots, n \text{ and } k=1,\ldots,2^m \}$  of $A^0_1$
up to $n^{\, \text{th}}$-level.    
This gives the corresponding Haar functions  
$\{ h_j \:  1\le j \le 2^n \}$.  
For each $ 1\le k \le 2^n$, 
we shall partition $A^n_k$ into 2 sets 
$A^{n+1}_{2k-1}$ and  $A^{n+1}_{2k}$ of equal measure 
(thus finding $h_{2^n+k}$) 
and find $x^*_{t_n} \in S(\X^*)$ and  
a sequence $\{ y^*_{n,i} \}_{i=1}^{p_n}$ in $S(\X^*)$ 
such that conditions (1), (2), and (3) hold.  

Find a finite set  $\{ y^*_{n,i} \}_{i=1}^{p_n}$ in $S(\X^*)$ 
that norms  $\text{sp} ( T h_j \:  1 \le j \le 2^n )$ 
within $\e_n$ . 
Let 
$$\eqalign{
  C^n_k  & \equiv   \cr
 & \{ f\in L_1 \: 
  |f| \le 1_{A^n_k} \text{, }\int_{A^n_k} f \, d\mu = 0 
  \text{ and } \int_{A^n_k}  (T^*y^*_{n,i}) f \, d\mu = 0 
  \text{ for }1 \le i \le p_n \}. 
}$$ 
Note that each $C^n_k$ is a convex weakly compact subset of $L_1$. 

Since $\{r_t\}$ tends weakly to $0$, for large $t$ there is a small
perturbation $\tilde r_t $ of $r_t$ so that $\tilde r_t 1_{A^n_k} $ 
 is in $C^n_k$ for each $k$.  To see this, put  
$$
  F = \text{sp}\left( \{ 1_{A^n_k} \}  \cup \{(T^*y^*_{n,i}) 1_{A^n_k}  
    \: k=1,\ldots, 2^n \text{ and } i=1,\ldots,p_n \} \right) \subset 
    L_\infty \simeq L_1^* \ . 
$$
Now
pick $t_n\equiv t$ so large that for $k=1,\ldots, 2^n$ and $i=1,\ldots,p_n$ 
\roster
\item"{(d)}" $\int_{A^n_k} r_t (T^*x^*_t) \, d\mu \ge 2\d\a_n$
\item"{(e)}" $\left|\int_\Omega r_t f \, d\mu \right| 
          \le  {\frac{\d}{2}} \a_n \n f \n $ \quad  for all $f$ in $F$
\endroster 
where $\a_n  = 2^{-n}~\mu(A^0_1) \equiv  \mu(A^n_k)$.  
Condition (d)
follows from the definition of $A_1^0$ and the weak-star convergence of 
$\{r_t(T^*x_t^*)\}$ to $h$ 
while condition (e) follows from (b) and the 
fact that $F$ is finite dimensional.

Thus the  $L_1$-distance from $r_t$ to 
$^\perp F \equiv 
\{ g\in L_1 \: \int_\Omega f g \, d\mu = 0 \text{ for each } f\in F \}$ 
is at 
most $\tfrac{\d~\a_n}{2}$.      
So there is $\tilde r_t \in  {^\perp F}$ such that 
$\n \tilde r_t - r_t \n_{L_1}$ is less than $\d\a_n $. 
Clearly $\tilde r_t 1_{A^n_k} \in C^n_k$ for each $k=1,\ldots, 2^n$. 

The functional $T^*x^*_t\in L^*_1$ attains its maximum on $C^n_k$ at an 
extreme point $u^n_k$ of $C^n_k$. 
By the lemma, $u^n_k = 1_{A^{n+1}_{2k-1}} - 1_{A^{n+1}_{2k}}$ 
for 2 disjoint sets $A^{n+1}_{2k-1}$ and $A^{n+1}_{2k}$ 
whose union is $A^n_k$.  
Furthermore, $A^{n+1}_{2k-1}$ and $A^{n+1}_{2k}$ are of 
equal measure since $\int_{A^n_k} u^n_k \, d\mu = 0$.  

Condition  (3) holds since for $i=1, \ldots, p_n $ and $k=1, \ldots,2^n$ 
$$
   y^*_{n,i}(Th_{2^n+k}) = 
     \a_n^{-1}~\int_{A^n_k} (T^*y^*_{n,i})u^n_k \, d\mu = 0 \ .
$$
Condition (1) follows from the observations that 
$$
x^*_{t_n} (T h_{2^n+k}) 
   =  \a_n^{-1} ~(T^* x^*_{t_n}) u^n_k
   \ge \a_n^{-1} ~(T^* x^*_{t_n}) (\tilde r_t 1_{A^n_k})
$$
and 
$$
| (T^* x^*_{t_n})(\tilde r_t 1_{A^n_k}) - (T^* x^*_{t_n})( r_t 1_{A^n_k}) | \le 
\n \tilde r_t - r_t\n_{L_1} < \d~\a_n  
$$
and 
$$
(T^* x^*_{t_n}) ( r_t 1_{A^n_k}) \ge  2\d\a_n \ .
\TagsOnRight \tag"\qed"\TagsOnLeft
$$
\enddemo

\demo{Proof of Lemma 2} 
Fix a function  $f$ of $C$ 
such that $|f| \neq 1_A $. 
Find  a positive $\alpha$ and a subset $B$ of $A$ with positive measure 
such that  $| f 1_B | < 1-\alpha$. 

Let $\tilde\Sg = B\cap\Sg$.  Consider the  measures 
$\lambda_i \: \tilde\Sg \to \Bbb R$  given  
by $\lambda_i (E) \equiv \int_E g_i \, d\mu$. 
Define the measure $\lambda \: \tilde\Sg \to \Bbb R^{n+1}$ by 
$$
  \lambda(E) = \left( \lambda_1 (E), \ldots, \lambda_n (E), \mu(E) \right) \ . 
$$
Liapounoff's Convexity Theorem gives a subset $B_1$  of $B$ satisfying 
\hfill\break
$\lambda(B_1) = \tfrac12\lambda(B) + \tfrac12 \lambda(\emptyset)$. 
Set $B_2 = B\setminus B_1$.  Note that 
$$
\lambda_i(B_1)  = \tfrac12\lambda_i(B)  = \lambda_i(B_2)  
\text{\qquad and\qquad}
\mu(B_1) = \tfrac12\mu(B) = \mu(B_2) \ 
$$
for $i=1,\ldots,n$.  Set 
$$
f_1 = f + \alpha\ (1_{B_1} - 1_{B_2}) 
\text{\qquad and \qquad}
f_2 = f + \alpha\ (1_{B_2} - 1_{B_1}) \ .
$$  
Clearly $f_1$ and $f_2$ are in $C$ 
and   $f = \tfrac{1}{2} f_1 + \tfrac{1}{2} f_2$. 
Thus $f$ is not an extreme point of $C$.\hfill\qed
\enddemo 

\heading{3. geometric view-point} \endheading \vskip 10pt

Consider a non-completely-continuous 
operator  $T\: L_1\to\X$ 
along with the corresponding 
Haar system $\{ h_j \}$ from Theorem~1. 
Let $\left\{ I^n_k  = [\tfrac{k-1}{2^n}, \tfrac{k}{2^n}) \right\}_{n,k}$ 
be the usual dyadic splitting of $[0,1]$ with corresponding 
Haar functions $\{ \tilde h_j \}_{j\ge 1}$. 
Consider the map $\tilde T \equiv T\circ S$ 
where $S\: L_1 \to L_1$ is the isometry 
that takes $\tilde h_j$ to $h_j$.  
Theorem~1 gives that there   
is a sequence $\{ x^*_n\}_{n\ge0}$ in $S(\X^*)$  
and a subspace $\X_0$ of $\X$  
such that
\roster 
\item $x^*_n ( \tilde T\tilde h_{2^n+k})  > \d$ for some $\d>0$ 
\item $\{ \text{sp} (\tilde T \tilde h_j \:  2^{n-1}<j\le 2^n ) \}_{n\ge0}$  
is a finite dimensional decomposition  of $\X_0$ with constant at most $1+\e$.  
\endroster
The next corollary follows from the observation that 
$\tilde T$ is not completely continuous and $\tilde T~L_1 \subset \X_0$. 
\proclaim{Corollary 3}
A Banach space failing the CCP has a 
subspace  with a finite dimensional decomposition (with constant 
arbitrarily close to 1) 
that fails the CCP.  
\endproclaim

A {\it tree} 
in  $\X$ is a system of the form    
\ $\{x^n_k  : n=0,1,\ldots \ ;\ k=1,\ldots,2^n \}$  \   
satisfying 
$$
x^{n}_k = \frac{x^{n+1}_{2k-1} + x^{n+1}_{2k}}{2}    \  .
$$
Associated to a tree is its difference system 
\ $\{d_j\}_{j\ge 1}$  \   where $d_1 = x^0_1$ and  
$$ 
 d_{2^n + k} \ =  \ \frac{x^{n+1}_{2k-1} \ - \ x^{n+1}_{2k}}{2}   \ . 
$$ 
There is a one-to-one correspondence between the bounded linear operators 
$T$ from $L_1$ into $\X$ and  bounded trees $\{x^n_k\}$ growing in $\X$. 
This correspondence is realized by  $T(\tilde h_j) = d_j$.  

A tree is a {\it $\d$-Rademacher tree}  if 
$\n \sum_{k=1}^{2^n} d_{2^n + k} \n ~\ge~2^n \d $.  
A tree  is a  {\it separated $\d$-tree}
if there exists a sequence 
$\{x^*_n \}_{n\ge 0}$ in $S(\X^*)$ such that 
$ x^*_n ( d_{2^n+k} ) \   >  \   \ \d  $. 
Clearly, a  separated $\d$-tree is also a $\d$-Rademacher tree.   
The operator corresponding to a $\d$-Rademacher tree
is not completely  
continuous since the image of the Rademacher functions 
stay large in norm.   
Thus if a  bounded $\d$-Rademacher tree (or  separated $\d$-tree) 
grows in $\X$, then  $\X$ fails the CCP.  

In any Banach space failing the CCP, 
a  bounded $\d$-Rademacher tree grows (see [G1] for a direct proof); 
in fact, even  
a bounded separated $\d$-tree grows  (see [G2] for an indirect proof). 
The proof of Theorem~1  is a {\it direct} proof that 
{\it if  $\X$ fails CCP then a bounded 
separated $\d$-tree, with a difference sequence 
naturally blocking into a finite dimensional decomposition, grows in $\X$.}  

As previously mentioned, we do not know whether a space failing CCP
must have a subspace with a basis which fails CCP.   
However, if  the space  also fails cotype   
(i.e. the space  contains $\ell_\infty^n$ uniformly for all $n$), 
then Theorem~1 can be modified  to show this is so.    

The main point is that if $\X$ fails cotype, 
$W$ is a finite dimensional subspace of $\X$, 
and $Z$ is a finite codimensional subspace of $\X$, then there is a
finite dimensional subspace $Y$ of $Z$  such that $W+Y$ has a basis with basis
constant less than, say, $10$. 
To see this, use the fact ([P], [JRZ]) that $W$ is 
$(1+\epsilon)$-complemented in a finite dimensional space
which has a basis with basis
constant less than $1+\epsilon$
and embed the complement to $W$ in that superspace
into $Z\cap {^\perp F}$, where $F$ is a finite subset of $\X^*$ which 
$(1+\epsilon)$-norms $W$.  
This is possible because finite codimensional subspaces of
$\X$ must contain $\ell_\infty^n$ uniformly for all $n$ and hence [J] contain
even $(1+\epsilon)$-isomorphs of $\ell_\infty^n$ for all $n$.

\proclaim{Corollary 4}
A Banach space failing the CCP  and failing cotype has a 
subspace  with a basis  
that fails the CCP.  
\endproclaim  

To see this,
it is enough by the argument for Corollary 3 to observe that
when $\X$  fails cotype  Theorem 1 can be improved by adding:
\roster
\item"{(D)}" a finite dimensional subspace $Y_n$ of $\X$ for each 
             $n\ge 0$ with $Y_0 = \emptyset$ , 
\endroster
changing (2) and (3) to:
\roster
\item"{(2$^\prime$)}" $\{ y^*_{n,i} \}_{i=1}^{p_n}$ norms  
      $\text{sp} (\cup_{k=0}^n Y_k \cup \{ T h_j \:  1 \le j \le 2^n \}) 
      $ within $\e_n$  
\item"{(3$^\prime$)}" $y^*_{n,i} (y) = 0 $ for
      $ y \in Y_{n+1} \cup \{ T h_{2^n + k}\}_{k=1}^{2^n} $ 
      and $1\le i \le p_n$,
\endroster
and adding:
\roster
\item"{(4)}" $\text{sp}\ (Y_{n+1} \cup \{ T h_{2^n + k}\}_{k=1}^{2^n}) $
      has a basis with basis constant less than $10$.
\endroster
Condition (2$^\prime$) is easily arranged. 
To achieve (D), (3$^\prime$), and (4), at the inductive step in the proof
of Theorem 1, after selecting $A_k^{n+1}$ (thereby defining $h_j$ for
$j=2^n+1,\dots,2^{n+1}$), choose   a finite dimensional 
 space $Y_{n+1} \subset \,
 ^\perp{\{}y_{n,i}^*\}_{i=1}^{p_n}$ so that 
$Y_{n+1} + \text{sp} \{ T h_{2^n + k}\}_{k=1}^{2^n}$ 
has a basis with basis constant less than $10$.

\widestnumber\no{[JRZ]Z}
\def\n #1{\no{[\bf #1]}}

\enddocument